\documentclass[a4paper]{amsart}
\usepackage[french]{babel}
\usepackage[latin1]{inputenc}
\usepackage{amssymb}
\usepackage{amsmath}
\usepackage{amscd}

\def\cal#1{\mathcal{#1}}
\def\AA{\mathbb{A}}
\def\PP{\mathbb{P}}

\def\CC{\mathbb{C}}

\def\NN{\mathbb{N}}
\def\ZZ{\mathbb{Z}}

\def\EE{\mathbb{E}}

\newtheorem{Theorem}{Theorem}[section]
\newtheorem{Corollary}[Theorem]{Corollary}
\newtheorem{Proposition}[Theorem]{Proposition}
\newtheorem{Lemma}[Theorem]{Lemma}
\newtheorem{Definition}[Theorem]{Definition}

\begin{document}

\title[Exceptional divisors which are not uniruled]{Exceptional divisors which are not uniruled
belong to the image of the Nash map}

\author[Monique Lejeune-Jalabert and Ana J. Reguera ]{Monique Lejeune-Jalabert and Ana J. Reguera}
\thanks{Partially supported by MTM2005-01518
} \maketitle

\centerline{\it Dedicated to H. Hironaka}

\vskip10mm

{\bf Abstract.} We prove that, if $X$ is a variety over an
uncountable algebraically closed field $k$ of characteristic zero,
then any irreducible exceptional divisor $E$ on a resolution of
singularities of $X$ which is not uniruled, belongs to the image
of the Nash map, i.e. corresponds to an irreducible component of
the space of arcs $X_\infty^\text{Sing}$ on $X$ centered in
$\text{Sing } X$. This reduces the Nash problem of arcs to
understanding which uniruled essential divisors are in the image
of the Nash map, more generally, how to determine the
uniruled essential divisors from the space of arcs. \\ \\

\section{Introduction}

In the midsixties, J. Nash ([Na]) initiated the study of the space
of arcs $X_\infty$ of a singular variety $X$ to understand what
the various resolutions of singularities of $X$ have in common.
His work was developed in the context of the proof of resolution
of singularities in characteristic zero by H. Hironaka ([Hi]).
From the existence of a resolution of singularities, Nash deduces
that the space of arcs $X_\infty^\text{Sing}$ on $X$ centered in
the singular locus $\text{Sing } X$ of $X$, has a finite number of
irreducible components. More precisely, he defines an injective
map ${\cal N}_X$, now called the {\it Nash map}, from the set of
irreducible components of $X_\infty^\text{Sing}$ which are not
contained in $(\text{Sing } X)_\infty$ to the set of {\it
essential divisors over} $X$, i.e. exceptional irreducible
divisors which appear up to birational equivalence on every
resolution of singularities of $X$. He asks whether this map is
surjective, or more generally, how complete is the description of
the essential
divisors by the image  of the Nash map. \\

In 1980, the first author ([Le]) proposed to approach the above
problem using arcs in the space of arcs $X_\infty$, or
equivalently, wedges. However, the ``curve selection lemma" does
not hold in $X_\infty$ because it is not a Noetherian space. This
obstacle was shortcut in 2006 by the second author ([Re2]), by
introducing the class of {\it generically stable irreducible
subsets of } $X_\infty$, and proving a curve selection lemma for
the corresponding {\it stable points of} $X_\infty$. The residue
field of such points being a transcendental extension of $k$ of
infinite transcendence degree, the stable points are very far from
being closed points. In [Re2] the image of the Nash map is
characterized in terms of a property of lifting wedges centered at
certain stable points to some resolution of singularities of $X$.
In this paper, given an essential divisor $\nu$ over $X$, we
introduce a property of {\it lifting wedges centered} at enough
closed points of a locally closed subset of $X_\infty$ associated
to $\nu$, which implies that $\nu$ belongs to the image of the
Nash map ${\cal
N}_X$ (see 2.1, def. 2.10 and cor. 2.15). \\

In 2003, S. Ishii and J. Kollar ([IK]) proved that the
surjectivity of the Nash map fails in general for $\dim X \geq 4$.
The fact that unirational does not imply rational for projective
varieties of dimension $\geq 3$ is crucial in the construction of
their example. The smooth cubic hypersurface in $\PP^4_\CC$ is the
first known example of a unirational variety which is not
rational, and as a consequence, it is uniruled but not
birationally ruled. In fact, if $k$ is algebraic closed, given a
smooth hypersurface $E$ in $\PP^{d}_k$ which is uniruled but not
birationally ruled, singular varieties $X$ of dimension $d \geq 4$
such that $E$ is an essential divisor over $X$ which is not in the
image of the Nash map are constructed in [IK]. \\

In section 3, we prove that, if the base field $k$ is
algebraically closed of characteristic zero and uncountable, then
any irreducible exceptional divisor $E$ on a resolution of
singularities $Y$ of $X$ which is not uniruled, is an essential
divisor over $X$ which belongs to the image of the Nash map (thm.
3.3). Using cor. 2.15, this follows from Luroth's theorem by
looking at the elimination of the points of indeterminacy of the
rational maps to $Y$ coming from the wedges. This reduces the Nash
problem for surfaces to decide which rational curves on the
minimal desingularization belong to the image of the Nash map.
Finally, in the appendix (prop. 4.2), we show that the
surjectivity of the Nash map for normal surface singularities over
the field $\CC$ of complex numbers would follow from proving that
every quasirational surface singularity over $\CC$ has a
resolution which enjoys the property of lifting wedges with
respect to each essential divisor. This makes more significant
studying
rational surface singularities ([Le], [Re1], [LR]).\\

{\it Acknowledgments: We would like to thank F. Loeser  and O.
Piltant
for their suggestions and encouragement.} \\

\section{Wedges and the image of the Nash map}

\begin{algo} \end{algo}
We begin this section by a brief introduction to the spaces of
arcs and wedges. For more details on arcs, see [DL], [EM], [IK],
[Vo].\\

Let $k$ be a field and let $X$ be a $k$-scheme. Given a field
extension $k \subseteq K$, a $K$-{\it arc} on $X$ is a
$k$-morphism $\text{Spec } K[[t]] \rightarrow X$. The $K$-arcs on
$X$ are the $K$-rational points of a $k$-scheme $X_\infty$ called
the {\it space of arcs} of $X$. More precisely, $X_\infty =
\lim_{\leftarrow} X_n$, where, for $n \in \NN$, $X_n$ is the
$k$-scheme of $n$-jets whose $K$-rational points are the
$k$-morphisms $\text{Spec } K[t]/ (t)^{n+1} \rightarrow X$. The
projective limit $\lim_{\leftarrow} X_n$ exists in the category of
$k$-schemes, because, for $n' \geq n$, the natural projections
$X_{n'} \rightarrow X_n$ are affine morphisms. For every
$k$-algebra $A$, we have a natural isomorphism
$$
Hom_k (\text{Spec } A, \ X_\infty) \cong Hom_k(\text{Spec }
A[[t]], \ X)  . \leqno(1)
$$

Given a $K$-arc $h: \text{Spec } K[[t]] \rightarrow X$, we will
call the image $h(O)$ in $X$ of the closed point $O$ of
$\text{Spec } K[[t]]$ the {\it center} of $h$. This is also the
image by the natural projection $j_0: X_\infty \rightarrow X$ of
the point of $X_\infty$ induced by the $K$-rational point of
$X_\infty$ corresponding to $h$ by (1). For $P \in X_\infty$ with
residue field $\kappa(P)$, we will denote by $h_P$ the
$\kappa(P)$-arc on $X$ corresponding by (1) to the
$\kappa(P)$-rational point of $X_\infty$ induced by $P$. We will
also call the center $j_0(P)$ of $h_P$ the center of $P$. Given $P
\in X_\infty$ as above with center $P_0$ on $X$, we will denote by
$\nu_P$ the {\it order function} $\text{ord}_t h_P^\sharp: {\cal
O}_{X,P_0} \rightarrow \NN \cup \{ \infty \}$ with $h_P^\sharp:
{\cal O}_{X,P_0} \rightarrow \kappa(P)[[t]]$ induced
by $h_P$.\\

A $K$-{\it wedge} on $X$ is a $k$-morphism $\text{Spec }
K[[\xi,t]] \rightarrow X$. The $K$-wedges on $X$ are the
$K$-rational points of a $k$-scheme $X_{\infty,\infty}$. We will
call $X_{\infty,\infty}$ the {\it space of wedges} of $X$. Since
by (1), for every $k$-algebra $A$, we have the natural
isomorphisms
$$
Hom_k (\text{Spec } A, \ (X_\infty)_\infty) \cong Hom_k(\text{Spec
} A[[\xi]], \ X_\infty) \cong Hom_k(\text{Spec } A[[\xi,t]], \ X)
 \leqno(2)
$$
the space of wedges $X_{\infty,\infty}$ and the space of arcs
$(X_\infty)_\infty$ of $X_\infty$ are naturally isomorphic
$k$-schemes, and we have a natural isomorphism
$$
Hom_k (\text{Spec } A, \ X_{\infty,\infty}) \cong Hom_k(\text{Spec
} A[[\xi,t]], \ X).
 \leqno(2')
$$
The $k$-scheme $X_{\infty, \infty}$ is the projective limit of the
$k$-schemes $\{X_{r,n}\}_{(r, n) \in \NN^2}$, where $X_{r,n}$ is
the $k$-scheme of $r$-jets of $X_n$, whose $K$-rational points are
the $k$-morphisms $\text{Spec } K[\xi,t]/ (\xi^{r+1},t^{n+1})
\rightarrow X$, with the natural affine transition morphisms
$X_{r',n'} \rightarrow X_{r,n}$, for $r' \geq r$ and $n' \geq
n$.\\

Given a $K$-wedge $\Phi: \text{Spec } K[[\xi,t]] \rightarrow X$,
we call the image in $X_\infty$ of the closed point (resp. generic
point) of $\text{Spec } K[[\xi]]$ by the $K$-arc $h_\Phi:
\text{Spec } K[[\xi]] \rightarrow X_\infty$ corresponding to
$\Phi$ by (2), the {\it special arc} (resp. {\it generic arc}) of
$\Phi$. With the above terminology, the special arc of $\Phi$ is
nothing but the center of the arc $h_\Phi$. This is also the image
by the natural projection $j_{0,\infty}: X_{\infty,\infty}=
\lim_{\leftarrow} X_{r,n} \rightarrow X_\infty=\lim_{\leftarrow}
X_{0,n}$ of the point of $X_{\infty, \infty}$ induced by the
$K$-rational point of $X_{\infty, \infty}$ corresponding to $\Phi$
by $(2')$. Thus, we will also say that $\Phi$ is {\it centered} at
$P \in X_\infty$ to mean that $P$ is the special arc of $\Phi$.\\

For ${\bf P} \in X_{\infty,\infty}$ with residue field
$\kappa({\bf P})$, we will denote by $\Phi_{\bf P}$ the
$\kappa({\bf P})$-wedge on $X$ corresponding by $(2')$ to the
$\kappa({\bf P})$-rational point of $X_{\infty,\infty}$ induced by
$\bf P$.\\

A $K$-wedge $\Phi$ on $X$ can also be viewed as a $K[[t]]$-point
of $X_\infty=\lim_{\leftarrow} X_{r,0}$. We will call the image in
$X_\infty$ of the closed point of $\text{Spec } K[[t]]$ the {\it
arc of centers} of $\Phi$. This is also the image by the natural
projection $j_{\infty,0}: X_{\infty,\infty}= \lim_{\leftarrow}
X_{r,n} \rightarrow X_\infty=\lim_{\leftarrow} X_{r,0}$ of the
point of $X_{\infty, \infty}$ induced by the $K$-rational point
of $X_{\infty, \infty}$ corresponding to $\Phi$ by $(2')$.\\

If $X$ is a variety over $k$, i.e. a reduced and irreducible
separated $k$-scheme of finite type, we will denote by
$X_\infty^\text{Sing}$ the closed set $j_0^{-1}(\text{Sing } X)$
of $X_\infty$ consisting of the arcs $P$ centered at some singular
point of $X$. Finally, we will denote by
$X_{\infty,\infty}^\text{Sing}$ the closed set $j_{\infty,0}^{-1}
(( \text{Sing } X)_\infty)$ consisting of the wedges $\bf P$ such
that the generic arc of $\Phi_{\bf P}$ belongs to
$X_\infty^\text{Sing}$.\\

Note that
$$
(\AA_k^{m})_\infty = \text{Spec } k[\{ {\underline {X}}_n \}_{n
\geq 0} ] \ \ \ \ \ (\AA_k^{m})_{\infty,\infty} = \text{Spec }
k[\{{\underline {X}}_{r,n}\}_{r,n \geq 0} ]
$$
where, ${\underline {X}}_n= ({X}_{1,n}, \ldots , {X}_{m,n})$ for
$n \geq 0$, and ${\underline {X}}_{r,n}= ({X}_{1,r,n}, \ldots ,
{X}_{m,r,n})$, for $r,n \geq 0$. Let us consider the morphism of
$k$-algebras ${\cal O}(\AA^{m}_k) \rightarrow {\cal
O}({(\AA^{m}_k)_\infty)[[t]]}$ induced in $(1)$ by the identity
map in $(\AA^{m}_k)_\infty$. For any $l \in {\cal O}(\AA_k^{m})$,
the image of $l$ in ${\cal O}((\AA^{m}_k)_\infty)[[t]]$ by the
previous morphism will be denoted by
$$
\sum_{n=0}^\infty L_n \ t^n \in {\cal O}((\AA_k^{m})_\infty) [[t]]
. \leqno (3)
$$
Analogously, for $l \in {\cal O}(\AA^{m}_k)$, the image of $l$ in
${\cal O}((\AA^{m}_k)_{\infty, \infty})[[\xi,t]]$ by the morphism
of $k$-algebras ${\cal O}(\AA^{m}_k) \rightarrow {\cal
O}({(\AA^{m}_k)_{\infty,\infty}})[[\xi,t]]$ induced in $(2')$ by
the identity map in $(\AA_k^m)_{\infty,\infty}$ will be denoted by
$$
\sum_{r,n \geq 0} L_{r,n} \ \xi^{r} \ t^{n} \in {\cal
O}((\AA_k^{m})_{\infty, \infty}) [[\xi, t]] . \leqno (4)
$$
If $X$ is a closed subscheme of $\AA^m_k$, and $I_X \subset k[x_1,
\ldots, x_m]$ is the ideal defining $X$ in $\AA^m_k$, then we have
$$
X_\infty= \text{Spec }  k [\{ {\underline {X}}_n, \}_{n \geq 0} ]
\ / \ (\{{F}_n\}_{n \geq 0, f \in I_X})
$$
$$
X_{\infty, \infty}= \text{Spec }  k [\{{\underline
{X}}_{r,n}\}_{r,n \geq 0}] \ / \ (\{{F}_{r,n}\}_{r,n \geq 0, f \in
I_X})
$$
We will use the same symbol to denote the element of $k[x_1,
\ldots, x_m]$, (resp.  ${\cal O}((\AA_k^{m})_\infty)$), (resp.
${\cal O}({(\AA_k^{m})_{\infty,\infty}})$) and its class in ${\cal
O}(X)$, (resp. ${\cal O}(X_\infty)$), (resp. ${\cal O}(X_{\infty, \infty})$).\\

\begin{algo}
\end{algo}
From now on, {\it $k$ will denote an uncountable field of
characteristic zero}, and $X$ a variety over $k$ of dimension $d$.
Let $p:Y \rightarrow X$ be a resolution of singularities of $X$,
i.e. a proper, birational morphism, with $Y$ nonsingular, such
that the induced morphism $Y \setminus p^{-1}(Sing \ X)
\rightarrow X \setminus Sing \ X$ is an isomorphism. Given an
irreducible component $E$ of the exceptional locus $p^{-1}(Sing \
X)$ of $p$ of codimension one, for short an {\it exceptional
divisor} from now on, we will denote by $Y^E_\infty$ the inverse
image of $E$ by the natural projection $j_0^Y: Y_\infty
\rightarrow Y$. Given a dense open subset $U$ of $E$ in the
nonsingular locus of $p^{-1}(Sing \ X)$, we will denote by
$Y_\infty^\dagger(U)$ the set of points $Q$ in $Y_\infty$ such
that $j_0^Y (Q) \in U$ and such that the corresponding arc $h_Q$
intersects $E$ transversally. Since $Y$ is nonsingular,
$Y^E_\infty$ and $Y_\infty^\dagger(U)$ are reduced and
irreducible. Moreover, $Y_\infty^\dagger(U)$ is open in
$Y^E_\infty$.\\

Let $p_\infty: Y_\infty \rightarrow X_\infty$ be the morphism
induced by $p$. The closure $N_E$ of $p_\infty(Y_\infty^E)$ is an
irreducible subset of $X_\infty^\text{Sing}$. Note that, if $P_E$
denotes the generic point of $N_E$, the order function $\nu_{P_E}$
coincides with the restriction of the divisorial valuation $\nu_E$
defined by $E$ to the local ring of $X$ at the generic point of
$p(E)$, and that, for every affine open subset $V$ of $X$ such
that $V \cap p(E) \not = \emptyset$, and for every $g \in {\cal
O}_X(V)$, we have
$$
\nu_{P_E}(g)= inf_{P \in N_E} \nu_P(g) .
$$
For $U$ as above, we set $N^\dagger(U) := p_\infty(
Y_\infty^\dagger(U))$. \\

\begin{Lemma}
Suppose in addition that $p$ is the blowing-up of a subscheme $D$
whose associated reduced scheme $D_\text{red}$ is $Sing \ X$, for
short an Hironaka resolution of singularities (see [Hi] I p. 132).
There exists a dense open subset $U$ of $E$ as above such that
$N^\dagger(U)$ is open in $N_E$; moreover, the morphism $p_\infty:
Y_\infty \rightarrow X_\infty$ induces an isomorphism of schemes
$Y_\infty^\dagger(U) \rightarrow N^\dagger(U)$.
\end{Lemma}

{\it Proof:} We may assume $X$ to be affine. Let $I=(g_0, \ldots ,
g_s) {\cal O}(X)$ be the ideal defining $D$ in $X$. We have $n_i:=
\nu_E(g_i) \geq 1$, for $0 \leq i \leq s$; for simplicity, suppose
that $n_0= inf \{ n_i \}_{0 \leq i \leq s}$.

Let $\Omega_0= \text{Spec } {\cal O}(X) [\frac {g_1}{g_0}, \ldots,
\frac {g_s}{g_0}]$ be the open subset of $Y$ on which $I {\cal
O}(\Omega_0)= (g_0) {\cal O}(\Omega_0)$, and let $U$ be the subset
of $\Omega_0 \cap E$ consisting of the nonsingular points of the
reduced exceptional locus of $p$. We will show that, for $P \in
N_E$, we have $P \in N^\dagger(U)$ if and only if
$\nu_P(g_0)=n_0$.

First note that, since $(div \ g_0)_U =n_0 \ E \cap U$, for $P \in
N^\dagger(U)$ we have $\nu_P(g_0)=n_0$. To prove the converse, we
first observe that, since for $P \in N_E$, we have $\nu_P(g_0)
\geq n_0$, the set $\{P \in N_E \ / \ \nu_P(g_0)=n_0\}$ is the
affine open set of $N_E$ where $(G_0)_{n_0} \neq 0$. Moreover,
since in addition, for $P \in N_E$, we have $\nu_P(g_i) \geq
\nu_{P_E}(g_i) =n_i \geq n_0$, for $0 \leq i \leq s$, the natural
morphism $\text{Spec } {\cal O}(N_E)_{(G_0)_{n_0}} [[t]]
\rightarrow X$ corresponding to the inclusion $N_E \setminus
V((G_0)_{n_0}) \hookrightarrow X_\infty$ by (1), lifts to
$\Omega_0 \subset Y$. As a consequence, the map sending $P \in N_E
\setminus V((G_0)_{n_0})$ to the center of the lifting
${\widetilde h}_P$ to $Y$ of the arc $h_P$, is a morphism of
schemes. Since this map sends $P_E$ to the generic point of $E$,
its image is contained in $E \cap \Omega_0$. We conclude that, for
any $P \in N_E \setminus V((G_0)_{n_0})$, the arc ${\widetilde
h}_P$ intersects transversally $E$ at a point in $\Omega_0$ which
is nonsingular in $p^{-1}(Sing \ X)$; so $N_E \setminus
V((G_0)_{n_0}) = N^\dagger(U)$.

Finally, the isomorphism $Y_\infty^\dagger(U) \rightarrow
N^\dagger(U)$ follows immediately, by (1) again, from the lifting
$\text{Spec } {\cal
O}(N_E)_{(G_0)_{n_0}} [[t]] \rightarrow Y$. \\

\begin{Remark}
\end{Remark}
In the previous lemma, we have used that, for a ring $R$, a series
$S$ in $R[[t]]$ is a unit if and only if its constant term $S(0)$
is invertible in $R$. In the next lemma, we will use the following
straightforward observation (which provides us with an algebraic
formula computing the coefficients in $t$ of the formal power
series
$y_i(t)=\frac{g_i(\underline{x}(t))}{g_0(\underline{x}(t))}$ from
the coefficients of $\underline{x}(t)$):

There exist polynomials $P_n \in$ $\ZZ[U_{n_0}, V_{n_0}, \ldots
,U_{n_0+n}, V_{n_0+n}]$, for $n \geq 0$, such that the equality
$$
\left( \sum_{n \geq 0} Y_n \ t^n \right) \ \left( \sum_{n \geq
n_0} U_n \ t^n \right) = \left( \sum_{n \geq n_0} V_n \ t^n
\right)
$$
in $k \left[\{Y_n\}_{n \geq 0}, \{U_n, V_n\}_{n \geq n_0} \right]
[[t]]$, is equivalent to
$$
Y_n \ {U_{n_0}^{n+1}= P_n(U_{n_0}, V_{n_0}, \ldots ,U_{n_0+n},
V_{n_0+n})},  \ \ \ \ \ \text{for } n \geq 0 . \leqno(5)
$$ \\


\begin{Lemma}
Let $p:Y \rightarrow X$ be an Hironaka resolution of singularities
of $X$. Let $k \subseteq K$ be a field extension, and let $\Phi$
be a $K$-wedge whose special arc is centered in $Sing \ X$, but
does not factor through $Sing \ X$.

If $\Phi$ does not lift to $Y$, then there exists a locally closed
subset $\Lambda$ of $X_{\infty, \infty}$ such that $\Phi$ is a
$K$-point of $\Lambda$, and such that, for any field extension $k
\subseteq L$, any $L$-point $\Psi$ of $\Lambda$ is a $L$-wedge on
$X$ which does not lift to $Y$.
\end{Lemma}

{\it Proof:} As in the proof of lemma 2.3, we may assume that $X$
is affine and that $p$ is the blowing-up of the ideal $I=(g_0,
\ldots , g_s) {\cal O}(X)$ with $V(I)=Sing \ X$.

First note that the special arc of any $L$-wedge $\Psi$ on $X$ is
centered in $Sing \ X$, but does not factor through $Sing \ X$, if
and only if the ideal $\Psi^\sharp(I)(0,t)$ in $L[[t]]$, image of
$\Psi^\sharp(I)$ by the map $L[[\xi,t]] \rightarrow L[[t]]$, $\xi
\mapsto 0, t \mapsto t$, is neither $(1)$ nor $(0)$, or
equivalently, $0 < ord_t \Psi^\sharp (I) (0,t) < \infty$.
Moreover, such a $\Psi$ lifts to $Y$ if and only if there exists
$i_0$, $0 \leq i_0 \leq s$, with $ord_t \Psi^\sharp(I) (0,t)=
ord_t \Psi^\sharp(g_{i_0}) (0,t)$, and such that the ideal
generated by $\Psi^\sharp(I)$ in $L[[\xi,t]]$ is the principal
ideal $(\Psi^\sharp \left(g_{i_0}) \right)$.

For simplicity suppose that $m_0:= ord_t \Phi^\sharp(g_0) (0,t)=
ord_t \Phi^\sharp(I) (0,t)$. Then $0 <m_0 < \infty$. The
$L$-wedges $\Psi$ such that $ord_t \Psi^\sharp(g_i) (0,t) \geq
m_0$ for $0 \leq i \leq s$, and $ord_t \Psi^\sharp(g_0)
(0,t)=m_0$, are the $L$-points of a locally closed subset
$\Lambda_{m_0}:=\{ (G_i)_{0,n}=0, \ 0 \leq i \leq s, \ 0 \leq n
<m_0, \ (G_0)_{0,m_0} \neq 0 \}$. Any $L$-point of $\Lambda_{m_0}$
is a $L$-wedge which lifts to $Y$ if and only if $\left(
\Psi^\sharp (I) \right)=\left( \Psi^\sharp(g_0) \right)$, and
$\Phi$ is a $K$-point of $\Lambda_{m_0}$.

Now suppose that $\Phi$ does not lift to $Y$, so there exists
$i_1$, $1 \leq i_1 \leq s$, such that $\Phi^\sharp(g_{i_1}) \notin
\left( \Phi^\sharp(g_0) \right)$. For simplicity, suppose that
$i_1=1$. Since $ord_t \Phi^\sharp(g_0)(0,t) = m_0 < \infty$, we
have that $n_0 := ord_t \Phi^\sharp(g_0) \leq m_0$, hence $e_0:=
ord_\xi \left( G_0^\Phi \right)_{n_0} < \infty$, where
$\Phi^\sharp(g_0)= \sum_{n \geq n_0} \left( G_0^\Phi \right)_n
t^n$ in $K[[\xi]][[t]]$. So either condition (i) or, in view of
(5), condition (ii) below holds:

\begin{itemize}

\item[(i)] $n_1:= ord_t \Phi^\sharp(g_1) < n_0$,

\item[(ii)] $ ord_t \Phi^\sharp(g_1) \geq n_0$ and there exists $n
\geq n_0$ such that
$$
e(n):= ord_\xi P_n \left( (G_0^\Phi)_{n_0}, (G_1^\Phi)_{n_0},
\ldots, (G_0^\Phi)_{n_0+n}, (G_1^\Phi)_{n_0+n} \right) < (n+1) e_0
,
$$
where $\Phi^\sharp(g_1)= \sum_{n \geq n_0} (G_1^\Phi)_{n} t^n$ in
$K[[\xi]][[t]]$.
\end{itemize}

In both cases, the $L$-wedges $\Psi$ such that $ord_t \Psi^\sharp
(g_0)=n_0$ and $ord_\xi \left( G_0^\Psi \right)_{n_0}=e_0$ are the
$L$-points of a locally closed subset $\Lambda^0_{n_0,e_0}:= \{
\left( G_0 \right)_{r,n}=0 \text{ for } 0 \leq n <n_0, \ \left(
G_0 \right)_{r,n_0}=0 \text{ for } 0 \leq r <e_0, \  \left( G_0
\right)_{e_0,n_0} \neq 0 \}$ in $X_{\infty \infty}$, and $\Phi$ is
a $K$-point of $\Lambda^0_{n_0,e_0}$.

In case (i), we have $e_1:=ord_\xi \left( G_1^\Phi \right)_{n_1} <
\infty$. Similarly, the $L$-wedges $\Psi$ such that $ord_t
\Psi^\sharp(g_1)=n_1$ and $ord_\xi \left( G_1^\Psi
\right)_{n_1}=e_1$ are the $L$-points of a locally closed subset
$\Lambda^1_{n_1, e_1}$ in $X_{\infty \infty}$, and $\Phi$ is a
$K$-point of $\Lambda_{n_1,e_1}^1$.

We conclude that $\Lambda:=\Lambda_{m_0} \cap \Lambda_{n_0,e_0}^0
\cap \Lambda_{n_1, e_1}^1$ is a locally closed subset in
$X_{\infty \infty}$, that $\Phi$ is a $K$-point of $\Lambda$ and
that any $L$-point $\Psi$ of $\Lambda$ is a $L$-wedge on $X$,
which does not lift to $Y$, since $\Psi^\sharp(g_1) /
\Psi^\sharp(g_0) \notin L[[\xi,t]]$.

In case (ii), let $n_2:=inf \{ n \in \NN \ / \ e(n) < (n+1) e_0
\}$. The $L$-wedges $\Psi$ such that

\begin{itemize}
\item[] $ord_t \Psi^\sharp(g_1) \geq n_0$,

\item[] $ord_\xi P_n \left( (G_0^\Psi)_{n_0}, (G_1^\Psi)_{n_0},
\ldots, (G_0^\Psi)_{n_0+n}, (G_1^\Psi)_{n_0+n} \right) \geq (n+1)
e_0$, for $0 \leq n <n_2$, \item[] $ord_\xi P_{n_2} \left(
(G_0^\Psi)_{n_0}, (G_1^\Psi)_{n_0}, \ldots, (G_0^\Psi)_{n_0+n_2},
(G_1^\Psi)_{n_0+n_2} \right) =e(n_2)$
\end{itemize}
are the $L$-points of a locally closed subset
$\Lambda^{0,1}_{n_0,e_0,e(n_2)}$ in $X_{\infty \infty}$ and $\Phi$
is a $K$-point of $\Lambda^{0,1}_{n_0,e_0,e(n_2)}$.

We conclude that $\Lambda := \Lambda_{m_0} \cap
\Lambda^{0}_{n_0,e_0} \cap \Lambda^{0,1}_{n_0,e_0,e(n_2)}$ is a
locally closed subset in $X_{\infty \infty}$, that $\Phi$ is a
$K$-point of $\Lambda$, and that any $L$-point $\Psi$ of $\Lambda$
is a $L$-wedge on $X$ which does not lift to $Y$, since again
$\Psi^\sharp (g_1) / \Psi^\sharp (g_0) \notin L[[\xi,t]]$.\\

\begin{Remark}
\end{Remark}
A similar statement holds for wedges which lift to $Y$, but we do
not need it in the sequel. \\

The following lemma
will be crucial in the proof of proposition 2.9. \\

\begin{Lemma}
Let $A$ be a noetherian integral ring and let $B$ be a countably
generated $A$-algebra. Let $\rho: Z=\text{Spec } B \rightarrow
S=\text{Spec }A$ be the induced morphism of affine schemes. If
$\Lambda$ is a locally closed subset of $Z$ such that the generic
point of $S$ lies in the image $\rho(\Lambda)$ of $\Lambda$ in
$S$, then there exists a countable family of dense open subsets
$\{U_n\}_{n \in \NN}$ of $S$ such that $\rho(\Lambda) \supseteq
\cap_{n \in \NN} U_n$.
\end{Lemma}

{\it Proof:} By assumption $\Lambda =U \cap F$ where $U$ is an
open set and $F$ is a closed set in $Z$. We have $F=\text{Spec }
B/I$ for an ideal $I$ in $B$, so $B/I$ is a countably generated
$A$-algebra. Now $\Lambda$ is an open set in $F$ and
$\rho'(\Lambda) = \rho(\Lambda)$ where $\rho': F \hookrightarrow Z
\rightarrow S$. This reduces to prove the lemma for an open subset
$\Lambda$ of $Z$.

The open sets of the form $D(g)$, for $g \in B$, complement of the
closed set $V((g))$, form a base of the Zariski topology of
$Z=\text{Spec } B$. So there exist $g_i \in B$, $i \in I$, such
that $\Lambda= \cup_{i \in I} D(g_i)$. We have $\rho(\Lambda) =
\cup_{i \in I} \rho(D(g_i))$, so there exists $i \in I$ such that
the generic point of $S$ lies in $\rho(D(g_i))$. This reduces to
prove the lemma for an open subset $\Lambda$ of $Z$ of the form
$D(g)$. But $D(g)= \text{Spec } B_g$ where $B_g$ is the localized
ring, which is isomorphic to $B[x]/(gx-1)$, hence a countably
generated $A$-algebra. This reduces to prove the lemma for
$\Lambda=Z$.

We have $B=\cup_{n \in \NN} B_n$ where $B_n$ is a finitely
generated $A$-algebra contained in $B$. Set $Z_n= \text{Spec }
B_n$ and let $\sigma_n : Z \rightarrow Z_n$ and $\rho_n: Z_n
\rightarrow S$ denote the canonical maps. For $n \in \NN$, we have
$\rho=\rho_n \circ \sigma_n$, so $\rho(Z) \subseteq \cap_n
\rho_n(Z_n)$. Now let $\wp$ be a point of $\text{Spec } A$ which
does not lie in $\rho(Z)$. The fiber of $\rho$ over this point is
empty, that is $\wp B_\wp=B_\wp$. So there exist an integer $e
\geq 1$, $a_i \in \wp$, $b_i \in B$, $1 \leq i \leq e$, and $s \in
A \setminus \wp$, such that $1 = \sum_{1 \leq i \leq e} a_i \ b_i
/s$ in $B_\wp$, i.e. there exists $s' \in A \setminus \wp$ such
that $s' \ (s- \sum_{1 \leq i \leq e} a_i \ b_i)=0$. Now, there
exists $n_0$ such that $b_i \in B_{n_0}$ for $1 \leq i \leq e$.
Hence $\wp (B_{n_0})_\wp =(B_{n_0})_\wp$, or equivalently $\wp$
does not lie in $\rho_{n_0}(Z_{n_0})$. We conclude that $\rho(Z)=
\cap_n \rho_n(Z_n)$. But, by Chevalley's theorem, each
$\rho_n(Z_n)$ is constructible in $S$. By assumption, the generic
point of $S$ lies in $\rho(Z)$, hence in $\rho_n(Z_n)$ for every
$n$, so $\rho_n(Z_n)$ contains a dense open subset $U_n$ of $S$
and $\rho(Z)$ contains $\cap_n U_n$.\\

\begin{Definition} ([Re3], definition 5.1)
Let $p:Y \rightarrow X$ be a resolution of singularities of $X$,
$E$ be an exceptional divisor on $Y$, and let $P_E$ be the generic
point of the irreducible closed subset $N_E$ of
$X_\infty^\text{Sing}$ defined in 2.2.

We say that $p$ satisfies the property of lifting wedges centered
at $P_E$ if, for any field extension $K$ of the residue field
$\kappa(P_E)$ of $P_E$ on $X_\infty$, any $K$-wedge on $X$ whose
special arc is $P_E$, lifts to $Y$.

We say that $p$ satisfies the property of lifting wedges in
$X_{\infty,\infty}^\text{Sing}$ centered at $P_E$ if, for any
field extension $K$ of $\kappa(P_E)$, any $K$-wedge on $X$ whose
special arc is $P_E$ and whose generic arc belongs to
$X_\infty^\text{Sing}$, lifts to $Y$.
\end{Definition}

Note that the property of lifting wedges centered at $P_E$ is
stronger that the property of lifting wedges in
$X_{\infty,\infty}^\text{Sing}$ centered at $P_E$.\\

\begin{Proposition}
Let $p:Y \rightarrow X$ be an Hironaka resolution of singularities
of $X$, let $E$ be an irreducible component of the exceptional
locus of $p$, and let $U$ be an open subset of $E$ as in lemma
2.3. Given  a point $Q_0 \in U$, we set
$Y_\infty^\dagger(Q_0)=Y_\infty^\dagger(U) \cap
(j_0^Y)^{-1}(Q_0)$, and
$N^\dagger(Q_0)=p_\infty(Y_\infty^\dagger(Q_0))$. Suppose that the
following property holds:
\begin{itemize}
\item[] ``For every countable family $\{U_n\}_{n \in \NN}$ of
nonempty open subsets of $U$, there exists a closed point $Q_0$ in
$\cap_n U_n$ such that, for any closed point $\bf P$ of
$X_{\infty, \infty}$ (resp. $X_{\infty, \infty}^\text{Sing}$) in
$j_{0, \infty}^{-1} (N^\dagger(Q_0))$ (i.e. such that the special
arc of the corresponding wedge $\Phi_{\bf P}$ on $X$ lies in
$N^\dagger (Q_0)$), the wedge $\Phi_{\bf P}$ lifts to $Y$."
\end{itemize}
Then, $p$ satisfies the property of lifting wedges centered at
$P_E$ (resp. $p$ satisfies the property of lifting wedges in
$X_{\infty,\infty}^\text{Sing}$ centered at $P_E$).
\end{Proposition}

{\it Proof:} Once again, we may assume $X$ to be affine. First
note that, since $p$ is an Hironaka resolution of singularities,
$E$ has codimension $1$ on $Y$. Suppose that $p$ does not satisfy
the property of lifting wedges centered at $P_E$ (resp. $p$ does
not satisfy the property of lifting wedges in
$X_{\infty,\infty}^\text{Sing}$ centered at $P_E$). Then, there
exist a field extension $K$ of $\kappa(P_E)$, and a $K$-wedge
$\Phi$ on $X$ (resp. a $K$-wedge $\Phi$ whose generic arc belongs
to  $X_\infty^{\text{Sing }}$) whose special arc is $P_E$ which
does not lift to $Y$. By lemma 2.5 there exists a locally closed
subset $\Lambda$ of $X_{\infty, \infty}$ such that $\Phi$ is a
$K$-point of $\Lambda$ and such that, for any field extension $k
\subset L$, any $L$-point of $\Lambda$ is a $L$-wedge on $X$ which
does not lift to $Y$. Let $\Omega_0$ be the affine open subset of
$Y$ and let $U$ be the open subset of $\Omega_0 \cap E$ defined in
the proof of lemma 2.3. Since $p_\infty$ induces an isomorphism of
schemes $Y^\dagger_\infty(U) \rightarrow N^\dagger (U)$, the map
$$
\rho: j_{0,\infty}^{-1}(N^\dagger(U)) \rightarrow N^\dagger(U)
\buildrel \sim \over{\longrightarrow} Y_\infty^\dagger(U)
\rightarrow U \subset \Omega_0 \cap E
$$
induced by the projections $j_{0,\infty}:X_{\infty, \infty}
\rightarrow X_\infty$, and $j_0^Y: Y_\infty \rightarrow Y$ is a
morphism of schemes. Now ${\cal O}(X_{\infty, \infty})$ is a
countably generated ${\cal O}(X_\infty)$-algebra via $j_{0,
\infty}$ and, setting $Y_\infty^{\Omega_0 \cap E}:=
(j_0^Y)^{-1}(\Omega_0 \cap E)$, \ ${\cal O}(Y_\infty^{\Omega_0
\cap E})$ is a countably generated ${\cal O} (\Omega_0 \cap
E)$-algebra via $j_0^Y$. Moreover $Y^\dagger_\infty
(U)=Y_\infty^{\Omega_0 \cap E} \setminus V((G_0)_{n_0})$ where, as
in lemma 2.3, $p$ is the blowing-up of the ideal $I$ of ${\cal
O}(X)$ and $I {\cal O}(\Omega_0)=(g_0) {\cal O}(\Omega_0)$. So,
the morphism $\rho$ enjoys the hypothesis of lemma 2.7. Let $
\Lambda':=\Lambda \cap j_{0,\infty}^{-1}(N^\dagger(U))$ (resp.
$:=\Lambda \cap j_{0,\infty}^{-1}(N^\dagger(U)) \cap
X_{\infty,\infty}^\text{Sing}$), which is a locally closed subset
of $X_{\infty, \infty}$. Moreover since $\Phi$ is a $K$-point of
$\Lambda'$ and the special arc of $\Phi$ is $P_E$, its image in
$\Lambda'$ is a point of $\Lambda'$ which is mapped to the generic
point of $\Omega_0 \cap E$ by $\rho$. Therefore, by lemma 2.7,
there exists a countable family of nonempty open subsets
$\{U_n\}_{n \in \NN}$ of $U$ such that $\cap_n U_n \subset
\rho(\Lambda')$. Since $k$ is uncountable, $\cap_n U_n$ has closed
points. For any such closed point $Q_0 \in \cap_n U_n$,
$\rho^{-1}(Q_0) \cap \Lambda'$ is a nonempty locally closed subset
of $X_{\infty, \infty}$. So the wedge $\Phi_{\bf P}$ corresponding
to any closed point ${\bf P} \in \rho^{-1}(Q_0) \cap \Lambda'$
does not lift to $Y$. This contradicts
the hypothesis.\\

The above proposition suggests to introduce:\\

\begin{Definition}
Let $p:Y \rightarrow X$ be a resolution of singularities of $X$,
and let $E$ be an exceptional divisor on $Y$. We say that $p$
satisfies the property of lifting wedges (resp. $p$ satisfies the
property of lifting wedges in $X_{\infty, \infty}^\text{Sing})$)
with respect to $E$, if the set of closed points $Q_0$ in $E$ such
that every wedge $\Phi_{\bf P}$ given by a closed point $\bf P$ of
$X_{\infty, \infty}$ (resp. $X_{\infty, \infty}^\text{Sing}$)
whose special arc lies in $N^\dagger(Q_0)$ lifts to $Y$, has a
nonempty intersection with $\cap_{n \in \NN} U_n$ for every family
of dense open subsets $\{U_n\}_{n \in \NN}$ of $E$.
\end{Definition}

\begin{algo}
\end{algo}
An {\it essential divisor} over $X$ is a divisorial valuation
$\nu$ of the function field $k(X)$ of $X$ centered in $\text{Sing
} X$, such that the center of $\nu$ on any resolution of
singularities $p: Y \rightarrow X$ is an irreducible component of
the exceptional locus of $p$ ([Na] p. 35, and [IK]). Given a
resolution of singularities $p: Y \rightarrow X$, and an
exceptional divisor $E$ on $Y$, if the divisorial valuation $\nu$
defined by $E$ is essential, then the point $P_{E}$ of $X_\infty$
defined in 2.2 only depends on $\nu$, not on $Y$. Let
$\{\nu_\alpha\}_{\alpha \in {\cal E}}$ be the set of essential
divisors over $X$, and for each $\alpha \in {\cal E}$, set
$P_\alpha:= P_{E_\alpha}$, where $E_\alpha$ is the center of
$\nu_\alpha$ in some divisorial resolution. Let
$N_\alpha:=\overline{\{P_\alpha\}}$, then
$\nu_\alpha:=\nu_{P_\alpha}$ only depends on $N_\alpha$.

Since $\text{char } k=0$, we have
$$
X_\infty^\text{Sing} = \bigcup_{\alpha \in {\cal E}} N_\alpha
$$
([Na], [IK], see also [Re2] prop. 2.2).\\

\begin{Definition}([Na]) We call the map
\[ {\cal N}_X : \ \ \{\text{irreducible components of
} X_\infty^{\text{Sing}}  \} \longrightarrow  \{ \text{essential
divisors over } X\}
\]
sending $N_\alpha \mapsto \nu_\alpha$ the ``Nash map". \\
\end{Definition}

The map ${\cal N}_X$ is injective but it need not be surjective if
$\dim X \geq 4$ ([IK]). Following [Le], we are looking for a
characterization of the image of ${\cal N}_X$ involving wedges.
The first approach appears in [Re2]:\\

\begin{Proposition} ([Re2], theorem 5.1)
Let $\nu_\alpha$ be an essential divisor over $X$, and let
$\kappa(P_\alpha)$ be the residue field of $P_\alpha$ in
$X_\infty$. The following conditions are equivalent:
\begin{itemize}
\item[(i)] $\nu_\alpha$ belongs to the image of the Nash map
${\cal N}_X$. \item[(ii)] For any resolution of singularities $p:
Y \rightarrow X$, $p$ satisfies the property of lifting wedges in
$X_{\infty,\infty}^\text{Sing}$ centered at $P_\alpha$, i.e. for
any field extension $K$ of $\kappa(P_\alpha)$, any $K$-wedge on
$X$ whose special arc is $P_\alpha$ and whose generic arc belongs
to $X_\infty^{Sing}$, lifts to $Y$. \item[(ii')] There exists a
resolution of singularities $p: Y \rightarrow X$ satisfying the
property of lifting wedges in
$X_{\infty,\infty}^\text{Sing}$ centered at $P_\alpha$.\\
\end{itemize}
\end{Proposition}

The previous result follows from the theory of {\it stable} points
of $X_\infty$ developed in [Re2] (see also [Re3]), whose key is to
understand the finiteness properties of the points $P_E$ defined
in 2.2. The
main result in [Re2] is the following:\\

\begin{algo}
\end{algo}
\textsc{Finiteness property of the stable point $P_E$}
\textit{([Re2], th. 4.1)}. The maximal ideal $P_E \ {\cal
O}_{{(X_\infty)_\text{red}},P_E}$ of the local ring ${\cal
O}_{{(X_\infty)_\text{red}},P_E}$ is finitely generated.
Therefore, $\widehat{{\cal O}_{{(X_\infty)_\text{red}},P_E}}$ is a
Noetherian complete local ring.\\

Taking into account that, for an uncountable algebraically closed
field $k$, the closed points of $X_{\infty, \infty}$ are the
$k$-rational points of $X_{\infty, \infty}$ (see [Is] prop. 2.10),
combining propositions 2.9 and 2.13, we conclude:\\

\begin{Corollary}
Let $k$ be an uncountable algebraically closed field of
characteristic zero. Let $\nu_\alpha$ be an essential divisor over
$X$, and let $p:Y \rightarrow X$ be an Hironaka resolution of
singularities.  Assume that the set of $k$-points $Q_0$ in
$E_\alpha$ such that every $k$-wedge $\Phi$ whose special arc lies
in $N^\dagger(Q_0)$ and whose generic arc belongs to
$X_\infty^\text{Sing}$ lifts to $Y$ has a nonempty intersection
with $\cap_{n \in \NN} U_n$ for every family of dense open subsets
$\{U_n\}_{n  \in \NN}$ of $E$. Then  $\nu_\alpha$ belongs to the
image of the Nash map ${\cal
N}_X$.\\\
\end{Corollary}

\section{Exceptional divisors which are not uniruled}

In this section $k$ is assumed to be an uncountable algebraically
closed field of characteristic zero.\\

\begin{Definition} Let $E$ be a variety over $k$  of
dimension $d$. We say that $E$ is uniruled if there exists a
variety $E'$ over $k$ of dimension $d-1$, and a dominant rational
map $\PP^1_k \times E' \cdots \rightarrow  E$.

\end{Definition}

\begin{Lemma} Let $E$ be an irreducible variety over $k$ which is
projective and nonsingular. Then $E$ is uniruled if and only if
there exists a countable family $\{U_n\}_{n \in \NN}$ of nonempty
open subsets of $E$ such that, for every $k$-point $Q_0$ in
$\cap_n U_n$, there is a $k$-rational curve in $E$ going through
$Q_0$.
\end{Lemma}

{\it Proof:} If $E$ is uniruled, there is a $k$-rational curve
through every $k$-point of $E$ ([Ko] chap.IV, cor. 1.4.4, [De]
chap. 4 rem. 4.2 (4)). Conversely, since $E$ is projective and
nonsingular, it is enough to show that $E$ contains a $k$-rational
curve $C$ over which the tangent bundle $T_E$ is generated by
global sections. Such a curve is said to be free ([Ko] chap. IV
th. 1.9, [De] chap. 4, cor. 4.11). Now, since $\text{char } k=0$,
by [Ko] chap. II th. 3.11, [De] chap. 4, prop. 4.14, there exists
a subset $E^\text{free}$ of $E$ and dense open subsets $V_n$, $n
\in \NN$, of $E$ such that $E^\text{free}= \cap_n V_n$ and such
that any $k$-rational curve $D$ with a nonempty intersection with
$E^\text{free}$ is free. Let $Q_0$ be a $k$-rational point of $E$
in $(\cap_n U_n ) \cap (\cap_n V_n)$. Such a point exists because
$k$ is uncountable and algebraically closed. By assumption, there
is a $k$-rational curve going through $Q_0$, hence it is free and
$E$ is uniruled.\\

\begin{Theorem} Let $X$ be a variety over an uncountable algebraically
closed field of characteristic zero $k$. Let  $p:Y \rightarrow X$
be a resolution of singularities of $X$, and let $E$ be an
exceptional divisor on $Y$.

If $E$ is not uniruled, then the divisorial valuation defined by
$E$ is an essential divisor which belongs to the image of the Nash
map.
\end{Theorem}

{\it Proof:} By [Ab1] prop. 4, if the valuation $\nu$ defined by
$E$ is not an essential divisor, then $E$ is birationally ruled,
i.e. $E$ is birationally isomorphic to $F \times \PP^1_k$ for some
$k$-variety $F$. In particular, $E$ is uniruled. So $\nu$ is an
essential divisor. Hence, we may assume that $p: Y \rightarrow X$
is an Hironaka resolution of singularities with nonsingular
exceptional divisors, since the center of $\nu$ on any such
resolution is a divisor birationally equivalent  to $E$, hence not
uniruled. Therefore, in view of prop. 2.13 (ii'), we only have to
prove that $p$ satisfies the property of lifting wedges in
$X_{\infty, \infty}^\text{Sing}$ centered at $P_E$. To do so, we
will apply cor. 2.15.

Let $U$ be an open subset of $E$ as in lemma 2.3, and let
$\{U_n\}_{n \in \NN}$ be a countable family of nonempty open
subsets of $U$. By lemma 3.2, since $E$ is nonsingular and not
uniruled, there exists a $k$-point $Q_0$ in $\cap_n U_n$, such
that no $k$-rational curve in $E$ goes through $Q_0$. Let $\Phi$
be a $k$-wedge on $X$ whose special arc is in $N^\dagger(Q_0)$. We
will prove that $\Phi$ lifts to $Y$.

Let $Z_0:=\text{Spec } k[[\xi,t]]$. The wedge $\Phi: Z_0
\rightarrow X$ does not factor through $\text{Sing } X$, because
the special arc of $\Phi$ does not. Hence there exists a finite
sequence of point blowing-ups $q:Z \rightarrow Z_0$ which
eliminates the points of indeterminacy of the rational map $Z_0
\dashrightarrow Y$ induced by $\Phi$,
so we have a commutative diagram of $k$-morphisms \\
$$
\begin{array}{ccc}
Z &  \buildrel  {\Phi'} \over{\longrightarrow}& Y\\
q \downarrow & &\downarrow  p\\
Z_0& \buildrel  {\Phi} \over{\longrightarrow}&X .\\ \\
\end{array}
$$
Let $C$ be the exceptional locus of $q$. Since $C$ is connected,
its image $\Phi'(C)$ by $\Phi': Z  \rightarrow  Y$ is a connected
subset of the exceptional locus of $p$ which contains $Q_0$. Since
$C$ is a finite union of $\PP^1_k$'s, by L\"{u}roth's theorem, our
hypothesis on $E$ implies that $\Phi'(C)$ is the point $Q_0$. Thus
$\Phi'$ factors through $\text{Spec } {\cal O}_{Y,Q_0}$, i.e. we
have
$$
\Phi':Z \longrightarrow \text{Spec } {\cal O}_{Y,Q_0} \ \subset \
Y .
$$
Indeed $\Phi'$ factors through every open subset $V$ of $Y$
containing $Q_0$, or equivalently $Z=\Phi'^{-1}(V)$. This is
because $Z \setminus \Phi'^{-1}(V)$ is a closed subset of $Z$
contained in $Z \setminus C$, hence $q$ induces an isomorphism
onto its image $F$ in $Z_0$. Thus $F$ is a closed subset of $Z_0$
which does not contain the closed point $O=q(C)$ of $Z_0$. This
implies that $F= \emptyset$, so $Z=\Phi'^{-1}(V)$.

Finally, since $Z_0$ is normal, we have that $q_\ast {\cal O}_Z
\cong {\cal O}_{Z_0}$ hence an isomorphism $k[[\xi,t]] \cong {\cal
O}(Z)$. So, for every $y \in {\cal O}_{Y,Q_0}$, its image in
${\cal O}(Z)$ lies in $k[[\xi,t]]$, i.e. there exists a morphism
$\widetilde{\Phi}:Z_0 \rightarrow \text{Spec } {\cal O}_{Y,Q_0}$
such that $\widetilde{\Phi} \circ q=\Phi'$. Since the morphism $q$
is birational, we get that $p \circ \widetilde{\Phi} =\Phi$. i.e.
$\Phi$ lifts
to $Y$.\\

\section{Appendix}

Combining Lur\"{o}th's and Grauert's contraction theorems, we show in
this appendix that a positive answer to the Nash problem for
normal surface singularities over $\CC$ would follow from proving
the property of lifting wedges with respect to the essential
divisors over those surface singularities over $\CC$ which are
quasirational.\\

\begin{Definition} ([Ab2])
Let $S$ be a surface over an algebraically closed field $k$ (i.e.
a variety of dimension $2$ over $k$), and let $P$ be a singular
point of $S$ at which $S$ is normal. We say that $S$ has a
quasirational singularity at $P$ if there exists a resolution of
singularities $p: Y \rightarrow (S,P):= \text{Spec } {\cal
O}_{S,P}$ such that the irreducible components of the exceptional
locus $p^{-1}(P)$ are rational curves.
\end{Definition}

Note that the same holds for every resolution of singularities of
$(S,P)$. Also note that, since the inverse image by the map
$\widehat{(S,P)}:=\text{Spec } \widehat{{\cal O}_{S,P}}
\rightarrow (S,P)$ of a resolution of singularites of $(S,P)$ is a
resolution of singularites of $\widehat{(S,P)}$ ([Li] lemma 16.1),
this is equivalent to saying that there exists a resolution of
singularities of $\widehat{(S,P)}$ (or equivalently for all) such
that the irreducible components of the exceptional locus are
rational curves. Here $\widehat{{\cal O}_{S,P}}$ is the completion
of
${\cal O}_{S,P}$ for the $M_{S,P}$-adic topology.\\

\begin{Proposition}
If every quasirational surface singularity over $\CC$ has a
resolution of singularities which satisfies the property of
lifting wedges with respect to each exceptional essential divisor,
then, for every normal surface $S$ over $\CC$, the Nash map ${\cal
N}_S$ is surjective.
\end{Proposition}

{\it Proof:} We may assume $S$ to be affine. Let $\nu$ be an
essential divisor over $S$, let $p:Y \rightarrow S$ be an Hironaka
resolution of singularities of $S$, and let $E$ be the exceptional
curve, center of $\nu$ on $Y$. If $E$ is not a rational curve, we
already know by theorem 3.3 that $\nu$ belongs to the image of the
Nash map. If not, let $P \in S$ denote the isolated singular point
of $S$ such that $p(E)=P$, let $\{E_\beta \}_{\beta \in B}$ be the
irreducible components of $p^{-1}(P)$ which are rational curves
and let $\EE$ be the connected component of $\cup_{\beta \in B}
E_\beta$ which contains $E$. The intersection matrix of the
irreducible curves in $\EE$ is negative definite. Hence by
Grauert's contraction theorem ([Gr] p. 367), there exists an
analytic normal surface $X$ and a proper holomorphic map
$p_1:Y^{\text{an}} \rightarrow X$ which contracts $\EE$ to a point
$Q \in X$ and induces a biholomorphic map from $Y^{\text{an}}
\setminus \EE$ to $X \setminus Q$. Here $Y^{\text{an}}$, below
$S^{\text{an}}$, denotes the associated complex analytic space.
The holomorphic map $X \setminus Q \cong Y^{\text{an}} \setminus
\EE \hookrightarrow Y^{\text{an}} \rightarrow S^{\text{an}}$
extends to a holomorphic map $p_2:X \rightarrow S^{\text{an}}$.
Indeed, since $S$ is normal, $S^{\text{an}}$ is again normal, so
for any open neighborhood $V$ of $P$ on $S^{\text{an}}$ (in the
complex topology) we have an isomorphism ${\cal
O}_{S^{\text{an}}}(V) \buildrel \sim \over{\longrightarrow} {\cal
O}_{Y^{\text{an}}}(p^{-1}(V))$. Now $V_1=p_1(p^{-1}(V))$ is an
open neighborhood of $Q$ on $X$ and $p_1^{-1}(V_1)=p^{-1}(V)$,
and, $X$ being normal, we also have an isomorphism ${\cal
O}_X(V_1) \buildrel \sim \over{\longrightarrow} {\cal
O}_{Y^{\text{an}}}(p_1^{-1}(V_1))$. The natural isomorphisms
${\cal O}_{S^{\text{an}}}(V) \cong {\cal O}_X(V_1)$ so obtained
give rise to a local morphism ${\cal O}_{S^{\text{an}},P}
\rightarrow {\cal O}_{X,Q}$, hence to a holomorphic map $p_2:X
\rightarrow S^{\text{an}}$ such that $p=p_2 \circ p_1$.\\

Let now $\Phi: Z_0:= \text{Spec } \CC[[\xi,t]] \rightarrow S$ be a
$\CC$-wedge on $S$, whose special arc lies in $N^\dagger(Q_0)$
with $Q_0$ a $\CC$-point of $E$ in the nonsingular locus of
$p^{-1}(P)_\text{red}$. We will first show that $\Phi$ lifts to
$X$. As in the proof of theorem 3.3, let $q:Z \rightarrow Z_0$ be
the finite sequence of point blowing-ups which eliminates the
points of indeterminacy of the rational map $Z_0 \cdots
\rightarrow Y$ induced by $\Phi$, let $C$ be the exceptional locus
of $q$ and let
\begin{center}
\begin{picture}(100,70)(0,0)
\put(0,55){\makebox(0,0)[br]{$Z$}}
\put(0,10){\makebox(0,0)[br]{$Z_0$}}
\put(95,10){\makebox(0,0)[br]{$S$ }}
\put(95,55){\makebox(0,0)[br]{$Y$ }}
\put(50,25){\makebox(0,0)[tr]{$\Phi$}}
\put(50,70){\makebox(0,0)[tr]{$\Phi'$}}
\put(-5,40){\makebox(0,0)[tr]{$q$}}
\put(100,40){\makebox(0,0)[tr]{$p$}} \put(30,55){\vector(1,0){30}}
\put(30,15){\vector(1,0){30}} \put(0, 45){\vector(0,-1){15}}
\put(85, 45){\vector(0,-1){15}}
\end{picture}
\end{center}
be the resulting commutative diagram. Here again $\Phi'(C)$ is a
connected subset of $p^{-1}(P)$ which has a nonempty intersection
with $E$, so by Lur\"{o}th's theorem, we have $\Phi'(C) \subset
\EE$.\\

Since $p_1^{-1}(Q)$ is a compact analytic subspace of
$p^{-1}(P)^\text{an}$ with reduced subspace $\EE^\text{an}$, and
$p^{-1}(P)$ is a projective scheme, by GAGA ([Se]) there exists an
ideal sheaf $\cal I$ on $Y$ such that $M_{X,Q} {\cal
O}_{Y^\text{an}}={\cal I} {\cal O}_{Y^\text{an}}$, and for every
$n \geq 1$, we have $\Gamma(Y^\text{an}, {\cal
O}_{Y^\text{an}}/{\cal I}^n{\cal O}_{Y^\text{an}}) = \Gamma (Y,
{\cal O}_Y / {\cal I}^n)$. Besides, by Grauert's comparison
theorem ([BS] chap. III sec.3), the natural map
$$
(p_1)_\ast({\cal O}_{Y^\text{an}})_Q \otimes_{{\cal O}_{X,Q}}
\widehat{{\cal O}_{X,Q}} \rightarrow \lim_\leftarrow
\Gamma(Y^\text{an}, {\cal O}_{Y^\text{an}}/{\cal I}^n{\cal
O}_{Y^\text{an}})
$$
is an isomorphism. But $X$ being normal, we have $(p_1)_\ast
({\cal O}_{Y^\text{an}})= {\cal O}_X$, so we get a natural
isomorphism
$$
\widehat{{\cal O}_{X,Q}} \buildrel \sim \over{\longrightarrow}
\lim_\leftarrow \Gamma (Y, {\cal O}_Y / {\cal I}^n)= \Gamma (Y,
\lim_\leftarrow {\cal O}_Y / {\cal I}^n)=: \Gamma
(\widehat{Y}_\EE, {\cal O}_{\widehat{Y}_\EE})
$$
Here $\widehat{{\cal O}_{X,Q}}$ is the completion of ${\cal
O}_{X,Q}$ for the $M_{X,Q}$-adic topology and $\widehat{Y}_\EE$ is
the formal completion of $Y$ along $\EE$.\\

Now, since $C \subset \Phi'^{-1}(\EE)$, there exists an effective
divisor $D$ whose support is $C$ and an ideal sheaf $\cal J$ on
$Z$ such that ${\cal I} {\cal O}_Z= {\cal J} {\cal O}_Z(-D)$. For
every $n \geq 1$, we thus get natural morphisms
$$
\Gamma (Y, {\cal O}_Y / {\cal I}^n) \rightarrow \Gamma (Z, {\cal
O}_Z / {\cal I}^n{\cal O}_Z) \rightarrow \Gamma (Z, {\cal O}_Z /
{\cal O}_Z(-n D))
$$
from which we derive a natural morphism $\widehat{{\cal O}_{X,Q}}
\rightarrow \lim_\leftarrow \Gamma(Z, {\cal O}_Z / {\cal
O}_Z(-nD))=: \Gamma(\widehat{Z}_C, {\cal O}_{\widehat{Z}_C})$
where $\widehat{Z}_C$ is the formal completion of $Z$ along $C$.
But $C$ is the support of both ${\cal O}_Z/M_{Z_0,O} {\cal O}_Z$
and ${\cal O}_Z / {\cal O}_Z(-D)$, so we have $\Gamma(
\widehat{Z}_C, {\cal O}_{\widehat{Z}_C}) = \lim_\leftarrow
\Gamma(Z,{\cal O}_Z / M_{Z_0,O}^n {\cal O}_Z)$, and
 since $q_\ast {\cal O}_Z={\cal O}_{Z_0}$, the natural map
$\widehat{{\cal O}_{Z_0,O}} \rightarrow \lim_\leftarrow \Gamma(Z,
{\cal O}_Z / M_{Z_0,O}^n {\cal O}_Z)$ is an isomorphism by the
theorem on formal functions ([Ha]). So we get a natural map
$\widehat{{\cal O}_{X,Q}} \rightarrow \widehat{{\cal O}_{Z_0,O}}$,
hence a $\CC$-wedge $\Psi: Z_0 \rightarrow \widehat{(X,Q)}:=
\text{Spec } \widehat{{\cal O}_{X,Q}}$.\\

Finally, we have $p_2 \circ \Psi=\Phi$, where here, for
simplicity, we have denoted $\Phi$ (resp. $p_2$) the morphism $Z_0
\rightarrow \widehat{(S,P)}:=\text{Spec } \widehat{{\cal
O}_{S,P}}$ (resp. $\widehat{(X,Q)} \rightarrow \widehat{(S,P)}$)
induced by $\Phi$ (resp. $p_2$). Indeed, once again by the theorem
on formal functions, the natural map $\widehat{{\cal O}_{S,P}}
\rightarrow \lim_\leftarrow \Gamma (Y, {\cal O}_Y / M_{S,P}^n
{\cal O}_Y ) =: \Gamma (\widehat{Y}_{p^{-1}(P)}, {\cal
O}_{\widehat{Y}, p^{-1}(P)})$ is an isomorphism, and the diagram
\begin{center}
\begin{picture}(100,70)(0,0)
\put(40,55){\makebox(0,0)[br]{$\widehat{Y}_\EE$}}
\put(0,10){\makebox(0,0)[br]{$\widehat{Z}_C$}}
\put(100,10){\makebox(0,0)[br]{$\widehat{Y}_{p^{-1}(P)}$ }}
\put(20,15){\vector(1,0){30}} \put(0,25){\vector(1,1){20}}
\put(40,45){\vector(1,-1){20}}
\end{picture}
\end{center}
is commutative.\\

Now, since the surface $X$ has an isolated singularity at $Q$,
there exists an (algebraic) surface $S_1$ over $\CC$ and a point
$P_1$ of $S_1$ such that $\widehat{{\cal O}_{S_1,P_1}}$ and
$\widehat{{\cal O}_{X,Q}}$ are $\CC$-isomorphic ([Ar] th. 3.8).
Moreover $S_1$ has a quasirational singularity at $P_1$. To prove
that, it is enough to check that $Q$ is a singular point of $X$,
and that there exists a resolution of singularities of
$\widehat{(X,Q)}$ such that the irreducible components of its
exceptional locus are rational curves. But, since $E$ is the
center on $Y$ of an essential divisor over $S$, its image on the
minimal resolution of singularities of $X$ is again a curve, hence
$Q$ is a singular
point of $X$.\\

Besides let $\widetilde Y$ denote the schematic closure of the
inverse image of $\widehat{(X,Q)} \setminus p_2^{-1}(P)$ in $Y
\times_S \widehat{(X,Q)}$ and let $\widehat{p_1}: \widetilde{Y}
\rightarrow \widehat{(X,Q)}$ denote the natural projection. The
map $\widehat{p_1}$ is a resolution of singularities of
$\widehat{(X,Q)}$ and its exceptional locus
$\widehat{p_1}^{-1}(Q)$ is $\CC$-isomorphic to $\EE$, hence its
irreducible components are rational curves. This also shows that
every resolution of singularities $\pi_1: Y_1 \rightarrow S_1$ has
an exceptional essential divisor $E_1$ birationally equivalent to
$E$. But, if $\pi_1$ satisfies the property of lifting wedges with
respect to $E_1$, the set of $\CC$-points $Q_0$ in the
intersection of a dense open set of $E$ isomorphic to a dense open
set of $E_1$ with the nonsingular locus of $p^{-1}(P)_\text{red}$,
such that every $\CC$-wedge $\Phi$ whose special arc lies in
$N^\dagger(Q_0)$ lifts to $Y$, has a nonempty intersection with
$\cap_{n \in \NN} U_n$ for every family of dense open subsets
$\{U_n\}_{n \in \NN}$ of $E$. By corollary 2.15, the essential
divisor $\nu=\nu_E$ over $S$ belongs to the image of the
Nash map ${\cal N}_S$.\\

\begin{Remark}
\end{Remark}
Note that the generic arc of the lifting $\Psi$ to
$\widehat{(X,Q)}$ of a $\CC$-wedge $\Phi$ on $\widehat{(S,P)}$
whose generic arc is centered at $P$ need not be centered at $Q$.
Therefore proving that every quasirational surface singularity
$(S_1, P_1)$ has a resolution of singularities which satisfies the
property of lifting wedges in ${S_1}_{\infty, \infty}^\text{Sing}$
with respect to each essential divisor would not be enough to
insure the surjectivity of the Nash map for every normal surface
over $\CC$.\\

\vskip5mm

\parindent=0mm
{\bf References.} \\
\begin{itemize}\parindent=10mm
\item[\hbox to\parindent{\enskip\mbox{[Ab1]}}] S.S. Abhyankar,
{\em On the valuations centered in a local domain}, Amer. J.
Math., 78, 321-348 (1956). \item[\hbox
to\parindent{\enskip\mbox{[Ab2]}}] S.S. Abhyankar, {\em
Quasirational singularities}, Amer. J. Math., 101, 267-300 (1979).

\item[\hbox to\parindent{\enskip\mbox{[Ar]}}] M. Artin, {\em
Algebraic approximation of structures over complete local rings},
Publ. Math. I.H.E.S. 36, 23-58 (1969).

\item[\hbox to\parindent{\enskip\mbox{[Bo]}}] C. Bouvier {\em
Diviseurs essentiels, composantes essentielles des vari\'{e}t\'{e}s
toriques singuli\`{e}res}, Duke Math. J. 91, no. 3, 609-620 (1998).

\item[\hbox to\parindent{\enskip\mbox{[BS]}}] C. B\v{a}nic\v{a},
O.St\v{a}n\v{a}\c{c}il\v{a}, {\em Algebraic methods in the global
theory of complex spaces}, John Wiley and Sons, London (1976).

\item[\hbox to\parindent{\enskip\mbox{[De]}}]O. Debarre, {\em
Higher-dimensional Algebraic Geometry}, Springer (2000).

\item[\hbox to\parindent{\enskip\mbox{[DL]}}] J. Denef, F. Loeser,
{\em Germs of arcs on singular algebraic varieties and motivic
integration}, Invent. math. 135, 201-232 (1999).

\item[\hbox to\parindent{\enskip\mbox{[EM]}}] L. Ein, M. Mustata,
{\em Jet schemes and singularities}, arXiv math AG/0612862 v2 to
appear in the Proceedings of the 2005 AMS Summer Research
Institute in Algebraic Geometry.

\item[\hbox to\parindent{\enskip\mbox{[Ha]}}] R. Hartshorne, {\em
Algebraic geometry}, Graduate Texts in Mathematics 52,
Springer-Verlag, New York (1977).

\item[\hbox to\parindent{\enskip\mbox{[Hi]}}] H. Hironaka, {\em
Resolution of singularities of an algebraic variety over a field
of characteristic zero}, Annals of Maths. 79; I, n.1, 109-203; II,
n.2, 205-326 (1964).

\item[\hbox to\parindent{\enskip\mbox{[GL]}}] G.
Gonzalez-Sprinberg, M. Lejeune-Jalabert, {\em Families of smooth
curves on surface singularities and wedges}, Ann. Polon. Math. 67
n.2, 179-190 (1997).

\item[\hbox to\parindent{\enskip\mbox{[Gr]}}] H. Grauert, {\em
\"{U}ber Modifikationen und exzeptionelle analytische Mengen}, Math.
Ann, 146, 331-368 (1962).

\item[\hbox to\parindent{\enskip\mbox{[Is]}}] S. Ishii, {\em The
arc space of a toric vriety}, J. Algebra 278, 666-683 (2004).

\item[\hbox to\parindent{\enskip\mbox{[IK]}}] S. Ishii, J. Koll\'{a}r,
{\em The Nash problem on arc families of singularities}, Duke
Math. J. 120, n.3, 601-620 (2003).

\item[\hbox to\parindent{\enskip\mbox{[Ko]}}] J. Koll\'{a}r, {\em
Rational curves on algebraic varieties}, Springer EMG 32 (1999).

\item[\hbox to\parindent{\enskip\mbox{[Le]}}] M. Lejeune-Jalabert,
{\em Arcs analytiques et r\'{e}solution minimale des singularit\'{e}s des
surfaces quasi-homog\`{e}nes}, LNM Springer 777, 303-336, (1980).

\item[\hbox to\parindent{\enskip\mbox{[LR]}}] M. Lejeune-Jalabert,
A. Reguera, {\em Arcs and wedges on sandwiched surface
singularities}, Amer. J. Math. 121, 1191-1213, (1999).

\item[\hbox to\parindent{\enskip\mbox{[Li]}}] J. Lipman, {\em
Rational singularities with applications to algebraic surfaces and
unique factorization}, Publ. Math. I.H.E.S. 36, 195-279 (1969).

\item[\hbox to\parindent{\enskip\mbox{[Na]}}] J. Nash, {\em Arc
structure of singularities}, Duke Math. J. 81, 207-212, (1995).

\item[\hbox to\parindent{\enskip\mbox{[Re1]}}] A.J. Reguera, {\em
Families of arcs on rational surface singularities}, Manuscripta
Math. 88, 321-333, (1995).

\item[\hbox to\parindent{\enskip\mbox{[Re2]}}] A.J. Reguera, {\em
A curve selection lemma in spaces of arcs and the image of the
Nash map}, Compositio Math. 142, 119-130, (2006).

\item[\hbox to\parindent{\enskip\mbox{[Re3]}}] A.J. Reguera, {\em
Towards the singular locus of the space of arcs}, To appear in
Amer. J. Math.

\item[\hbox to\parindent{\enskip\mbox{[Se]}}] J.P. Serre, {\em
G\'{e}om\'{e}trie alg\'{e}brique et g\'{e}om\'{e}trie analytique}, Ann. Inst. Fourier
Grenoble 6, 1-42 (1956).

\item[\hbox to\parindent{\enskip\mbox{[Vo]}}] P. Vojta, {\em Jets
via Hasse-Schmidt derivations}, arXiv math.AG/0407113 v2.
\par\end{itemize}

\ \\ \\
\noindent Monique Lejeune-Jalabert,\\ CNRS, LMV, UMR8100,
CNRS-UVSQ, Universit\'{e} de Versailles, Saint-Quentin, B\^{a}timent
Fermat, 45 Avenue des Etats-Unis,
F-78035 Versailles Cedex, France.\\
E-mail: lejeune@math.uvsq.fr
\ \\ \\
\noindent
Ana J. Reguera, \\ Dep. de \'Algebra, Geometr\'ia y Topolog\'ia,  Universidad de Valladolid,\\
Prado de la Magdalena s/n,
47005 Valladolid, Spain. \\
E-mail: areguera@agt.uva.es

\end{document}